\documentclass{gtart}
\usepackage{curvesls}
\usepackage{epic}

\theoremstyle{plain}
\newtheorem{thm}{Theorem}[section]
\newtheorem{lem}[thm]{Lemma}

\theoremstyle{definition}
\newtheorem{defn}[thm]{Definition}

\newcommand{\Z}{{\bf{Z}}}
\newcommand{\C}{{\bf{C}}}
\newcommand{\Q}{{\bf{Q}}}
\newcommand{\cB}{{\mathcal{B}}}
\newcommand{\cH}{{\mathcal{H}}}

\newcommand{\TL}{{\mathrm{TL}}}

\newcommand{\half}{{\frac{1}{2}}}

\newcommand{\nchoosetwo}{n \choose 2}

\newcommand{\sijclose}
{\begin{picture}(110,120)
  \thinlines
  \put(0,0){\line(0,1){110}}
  \put(0,0){\line(1,0){110}}
  \put(110,110){\line(0,-1){110}}
  \put(110,110){\line(-1,0){110}}

  \thicklines

  \multiput(15,0)(20,0){2}{\oval(10,15)[t]}

  \multiput(0,0)(20,0){2}
  {
    \qbezier(10,110)(10,75)(30,55)
    \qbezier(30,55)(50,35)(50,0)
  }

  \multiput(35,55)(5,0){3}{\circle*{1}}

  \put(55,110){\oval(10,15)[b]}
  \put(55,110){\oval(30,35)[b]}

  \multiput(80,0)(20,0){2}{\line(0,1){110}}
  \multiput(85,55)(5,0){3}{\circle*{1}}

  \put( 5,110){\makebox(10,10){$1$}}
  \put(45,110){\makebox(10,10){$i$}}
  \put(65,110){\makebox(10,10){$j$}}
  \put(95,110){\makebox(10,10){$n$}}

\end{picture}
}

\newcommand{\sijfar}
{\begin{picture}(110,120)
  \thinlines
  \put(0,0){\line(0,1){110}}
  \put(0,0){\line(1,0){110}}
  \put(110,110){\line(0,-1){110}}
  \put(110,110){\line(-1,0){110}}

  \thicklines

  \multiput(15,0)(20,0){2}{\oval(10,15)[t]}

  \multiput(0,0)(10,0){2}
  {
    \qbezier(10,110)(10,75)(30,55)
    \qbezier(30,55)(50,35)(50,0)
  }
  \multiput(0,0)(10,0){2}
  {
    \qbezier(50,110)(50,65)(60,55)
    \qbezier(60,55)(70,45)(70,0)
  }
  \multiput(90,0)(10,0){2}{\line(0,1){110}}

  \multiput(35,110)(40,0){2}{\oval(10,15)[b]}

  \multiput(32.5,55)(2.5,0){3}{\circle*{1}}
  \multiput(62.5,55)(2.5,0){3}{\circle*{1}}
  \multiput(92.5,55)(2.5,0){3}{\circle*{1}}

  \put( 5,110){\makebox(10,10){$1$}}
  \put(25,110){\makebox(10,10){$i$}}
  \put(75,110){\makebox(10,10){$j$}}
  \put(95,110){\makebox(10,10){$n$}}
\end{picture}
}

\newcommand{\figeight}
{
  \qbezier(-6,  0)(-6, 10)( 0, 10)
  \qbezier( 0, 10)( 6, 10)(10,  0)
  \qbezier(10,  0)(14,-10)(20,-10)
  \qbezier(20,-10)(26,-10)(26,  0)
  \qbezier(26,  0)(26, 10)(20, 10)
  \qbezier(20, 10)(14, 10)(10,  0)
  \qbezier(10,  0)( 6,-10)( 0,-10)
  \qbezier( 0,-10)(-6,-10)(-6,  0)
}

\newcommand{\oooo}
{\raisebox{-48pt}
  {
  \begin{picture}(100,100)
  \put(50,50){\bigcircle{100}}
  \multiput(20,50)(20,0){4}{\circle*{4}}
  \multiput(20,50)(40,0){2} {\figeight}
  \end{picture}
  }
}
\newcommand{\oooosquare}
{\raisebox{-48pt}
  {
  \begin{picture}(100,100)
  \put(50,50){\bigcircle{100}}
  \multiput(20,50)(20,0){4}{\circle*{4}}
  \multiput(20,50)(40,0){2}{\line(1,0){20}}
  \end{picture}
  }
}

\newcommand{\ooo}
{\raisebox{-38pt}
  {
  \begin{picture}(80,80)
  \put(40,40){\bigcircle{80}}
  \multiput(20,40)(20,0){3}{\circle*{4}}
  \multiput(18,40)(24,0){2} {\figeight}
  \end{picture}
  }
}
\newcommand{\ooosquare}
{\raisebox{-38pt}
  {
  \begin{picture}(80,80)
  \put(40,40){\bigcircle{80}}
  \multiput(20,40)(20,0){3}{\circle*{4}}
  \put(20,40){\line(1,0){40}}
  \end{picture}
  }
}
\newcommand{\ooodisk}
{\raisebox{-38pt}
  {
  \begin{picture}(60,60)
  \put(30,30){\bigcircle{60}}
  \put(30,30){\circle*{4}}
  \put(20,1.72){\vector(1,1){20}}
  \qbezier(40,21.72)(48.38, 30)(40, 38.38)
  \put(40,38.38){\line(-1,1){20}}
  \put(40,58.38){\vector(-1,-1){20}}
  \qbezier(20, 38.38)(11.62, 30)(20, 21.72)
  \put(20,21.72){\line(1,-1){20}}
  \end{picture}
  }
}

\newcommand{\oo}
{\raisebox{-28pt}
  {
  \begin{picture}(60,60)
  \put(30,30){\bigcircle{60}}
  \multiput(20,30)(20,0){2}{\circle*{4}}
  \multiput(20,25)(0,10){2} {\figeight}
  \end{picture}
  }
}

\newcommand{\ootriangle}
{\raisebox{-28pt}
  {
  \begin{picture}(60,60)
  \put(30,30){\bigcircle{60}}
  \multiput(20,30)(20,0){2}{\circle*{4}}
  \put(20,30){\line(1,0){20}}
  \end{picture}
  }
}

\newcommand{\mm}
{\raisebox{-28pt}
  {
  \begin{picture}(60,60)
  \put(30,30){\bigcircle{60}}
  \multiput(10,30)(10,0){5}{\circle*{2}}

  \put(7,10.7){\line(0,1){20.3}}
  \put(10,30){\oval(6,10)[t]}
  \put(13,30){\line(0,-1){15}}
  \put(15,15){\oval(4,4)[b]}
  \put(17,15){\line(0,1){15}}
  \put(20,30){\oval(6,10)[t]}
  \put(23,30){\line(0,-1){29.2}}

  \put(53,10.7){\line(0,1){20.3}}
  \put(50,30){\oval(6,10)[t]}
  \put(47,30){\line(0,-1){15}}
  \put(45,15){\oval(4,4)[b]}
  \put(43,15){\line(0,1){15}}
  \put(40,30){\oval(6,10)[t]}
  \put(37,30){\line(0,-1){29.2}}
  \end{picture}

  }
}
\newcommand{\chopstickstwofour}
{\raisebox{-28pt}
  {
  \begin{picture}(60,60)
  \put(30,30){\bigcircle{60}}
  \multiput(10,30)(10,0){5}{\circle*{2}}
  \drawline(25,0.42)(20,30)
  \drawline(35,0.42)(40,30)
  \end{picture}
  }
}
\newcommand{\chopsticksonefour}
{\raisebox{-28pt}
  {
  \begin{picture}(60,60)
  \put(30,30){\bigcircle{60}}
  \multiput(10,30)(10,0){5}{\circle*{2}}
  \drawline(25,0.42)(10,30)
  \drawline(35,0.42)(40,30)
  \end{picture}
  }
}

\newcommand{\unionjack}
{
  \begin{picture}(100,100)
  \multiput(0,0)(50,0){3}{\line(0,1){100}}
  \multiput(0,0)(0,50){3}{\line(1,0){100}}
  \put(0,0){\line(1,1){100}}
  \put(100,0){\line(-1,1){100}}

  \put(10,30){\makebox(10,10){$1$}}
  \put(10,60){\makebox(10,10){$-q$}}
  \put(30,10){\makebox(10,10){$-t$}}
  \put(30,80){\makebox(10,10){$qt$}}
  \put(60,10){\makebox(10,10){$qt$}}
  \put(60,80){\makebox(10,10){$-q^2t$}}
  \put(80,30){\makebox(10,10){$-qt^2$}}
  \put(80,60){\makebox(10,10){$q^2t^2$}}
  \end{picture}
}

\newcommand{\xonen}
{\raisebox{-28pt}
  {
  \begin{picture}(60,60)
  \put(30,30){\bigcircle{60}}
  \multiput(12,30)(12,0){4}{\circle*{4}}

  \put(18,2.5){\line(0,1){55}}
  \put(18,30){\vector(0,1){0}}

  \qbezier(42,2.5)(30,14.5)(30,30)
  \put(36,30){\oval(12,18)[t]}
  \put(42,30){\vector(0,-1){0}}
  \put(48,30){\oval(12,18)[b]}
  \qbezier(54,30)(54,45.5)(42,57.5)
  \end{picture}
  }
}

\newcommand{\yone}
{
  \begin{picture}(100,100)
  \put(50,50){\bigcircle{100}}
  \multiput(20,50)(20,0){4}{\circle*{4}}
  \put(30,4.2){\line(0,1){91.6}}
  \put(60,50){\figeight}
  \end{picture}
}
\begin{document}

\title{The Lawrence-Krammer representation}
\author{Stephen J. Bigelow}
\address{Stephen J. Bigelow \\
         Department of Mathematics and Statistics \\
         University of Melbourne \\
         Victoria 3010 \\
         Australia}
\email{bigelow@unimelb.edu.au}

\begin{abstract}
The Lawrence-Krammer representation of the braid groups
recently came to prominence when it was shown to be faithful 
by myself and Krammer.
It is an action of the braid group
on a certain homology module $H_2(\tilde{C})$
over the ring of Laurent polynomials in $q$ and $t$.
In this paper we describe some surfaces in $\tilde{C}$
representing elements of homology.
We use these to give a new proof that $H_2(\tilde{C})$ is a free module.
We also show that the $(n-2,2)$ representation of the Temperley-Lieb algebra
is the image of a map to relative homology at $t=-q^{-1}$,
clarifying work of Lawrence.
\end{abstract}


\makeshorttitle

\section{Introduction}

The Lawrence-Krammer representation
is the action of the braid group $B_n$
on a certain homology module $H_2(\tilde{C})$ over the ring
$\Lambda = \Z[q^{\pm 1},t^{\pm 1}]$.
It was introduced by Lawrence \cite{rL90},
and recently came to prominence when 
it was shown to be faithful for $n=4$ \cite{dK00},
and then for all $n$ \cite{sB01}, \cite{dK02}.
Since then, a number of papers have appeared
which closely examine certain aspects of the representation.
See \cite{PP02}, \cite{wS02}, and \cite{rB02}.
We now continue in this tradition,
with the specific goal of understanding
the connection with the Temperley-Lieb algebra.

This paper was partly motivated by 
an attempt to clarify two points from \cite{sB01}.
First, the pairing between a ``noodle'' and a ``fork''
involved an algebraic intersection number
between two non-compact surfaces in $\tilde{C}$.
Such a thing is not necessarily well-defined,
so I gave an indirect proof of the existence of
a certain closed surface corresponding to a fork.
I now have an explicit description of this surface,
which will be given in Section \ref{sec:genus_three}.
In fact it is possible to define the pairing 
without reference to this surface
by using results from \cite[Appendix E]{aK96}.
(Thanks to Won Taek Song,
whose paper \cite{wS02} drew my attention to this.)

Second, to compute matrices for the representation,
I tensored $H_2(\tilde{C})$ with a field containing $\Lambda$.
The resulting vector space contains $H_2(\tilde{C})$,
but strictly speaking, the action of $B_n$ on this vector space
should not be called the Lawrence-Krammer representation.
In Section \ref{sec:basis} we give a new proof
that $H_2(\tilde{C})$ is a free $\Lambda$-module.
This is originally due to Paoluzzi and Paris \cite{PP02},
but our proof uses an explicit description of
surfaces representing elements of a free basis of $H_2(\tilde{C})$.
The correct matrices for the Lawrence-Krammer representation
could be computed from this,
but we do not do so since they are quite complicated.

By addressing these two issues from \cite{sB01}
we also shed new light on the work of Lawrence \cite{rL90}.
There, the Lawrence-Krammer representation
was used to give a topological interpretation
to a representation of the Temperley-Lieb algebra.
The idea is to specialise $t$ to the value $-q^{-1}$,
at which point the Lawrence-Krammer representation becomes reducible
and the desired Temperley-Lieb representation appears as a quotient.
Somewhat complicated methods are used in \cite{rL90}
to define the required quotient.
In Section \ref{sec:tl} we show that it is simply the image
of a map to a certain relative homology module.
The significance of the value $t = -q^{-1}$
can then be seen to follow from
the existence of certain kind of surface in $\tilde{C}$
that maps to a multiple of $1+qt$ in relative homology.

Throughout this paper $n$ is a positive integer,
$D$ is the unit disk centred at the origin in the complex plane,
$-1 < p_1 < \dots < p_n < 1$ are real numbers,
and
$D_n = D \setminus \{p_1,\dots,p_n\}$ is the $n$-times punctured disk.
The braid group $B_n$ is the mapping class group of $D_n$.
We also assume familiarity with 
the presentation of $B_n$ using Artin generators $\sigma_i$,
and the interpretation of $B_n$
as the fundamental group of a configuration space in the plane.

This paper is based on a talk I gave at
the 2001 Georgia International Topology Conference.
The research was supported by the Australian Research Council.

\section{Definitions}

\subsection{The Lawrence-Krammer representation}

We briefly review the definition of the Lawrence-Krammer representation
as given in \cite{sB01}.
Let $C$ be the space of unordered pairs of distinct points in $D_n$.
Let $c_0 = \{d_1,d_2\}$ be a basepoint in $C$,
where $d_1$ and $d_2$ are distinct points on the boundary of the disk.

We define a homomorphism
$$\Phi \co \pi_1(C,c_0) \to \langle q \rangle \oplus \langle t \rangle$$
as follows.
Suppose $\alpha \co I \to C$ is a closed loop in $C$
representing an element of $\pi_1(C,c_0)$.
By ignoring the puncture points we can consider
$\alpha$ as a loop in the space of unordered pairs of points in the disk,
and hence as a braid in $B_2$.
Let $b$ be the exponent of this braid
in the Artin generator $\sigma_1$.
Similarly, the map
$$s \mapsto \{p_1,\dots,p_{n}\} \cup \alpha(s)$$
determines a braid in $B_{n+2}$.
Let $b'$ be the the exponent sum of this braid
in the Artin generators $\sigma_i$ of $B_{n+2}$.
Note that $b$ and $b'$ have the same parity.
Let $a = \half(b'-b)$.
We define
$$\Phi(\alpha) = q^a t^b.$$

Let $\tilde{C}$ be the connected covering space of $C$ such that
$\pi_1(\tilde{C}) = \ker(\Phi)$.
Fix a choice of $\tilde{c}_0$ in the fibre over $c_0$.
The homology group $H_2(\tilde{C})$ admits a $\Lambda$-module structure.
Suppose $f$ is a homeomorphism from $D_n$ to itself, representing an
element of $B_n$.
The induced map $f_*$ from $\pi_1(C,c_0)$ to itself has the property
$$\Phi \circ f_* = \Phi.$$
It follows that $f$ has a unique lift 
$$\tilde{f} \co (\tilde{C},\tilde{c}_0) \to (\tilde{C},\tilde{c}_0).$$
Further, the induced map
$$\tilde{f}_* \co H_2(\tilde{C}) \to H_2(\tilde{C})$$
is a $\Lambda$-module automorphism.
We define the Lawrence-Krammer representation to be
this action of $B_n$ on $H_2(\tilde{C})$.

\subsection{An intersection pairing}
\label{sec:pairing}

We now introduce some relative homology modules
and an intersection pairing.
For $\epsilon > 0$, let $\nu_\epsilon$ be the set of points
$\{x,y\} \in C$ such that
either $x$ and $y$ are within distance $\epsilon$ of each other, 
or at least one of them is within distance $\epsilon$ 
of a puncture point.
Let $\tilde{\nu}_\epsilon$ be the preimage of $\nu_\epsilon$ in $\tilde{C}$.
The relative homology modules $H_2(\tilde{C},\tilde{\nu}_\epsilon)$ 
are nested by inclusion.
Let
$$H_2(\tilde{C},\tilde{\nu}) = \lim_{\epsilon \to 0} 
  H_2(\tilde{C},\tilde{\nu}_\epsilon),$$
and
$$H_2(\tilde{C},\partial\tilde{C} \cup \tilde{\nu}) = \lim_{\epsilon \to 0} 
  H_2(\tilde{C},\partial\tilde{C} \cup \tilde{\nu}_\epsilon).$$
The braid group $B_n$ acts on these,
and on $H_2(\tilde{C},\partial\tilde{C})$,
by $\Lambda$-module automorphisms.

For $x \in H_2(\tilde{C})$ and 
$y \in H_2(\tilde{C},\partial\tilde{C} \cup \tilde{\nu})$
let $(x \cdot y) \in \Z$ denote the standard intersection number.
We define an intersection pairing
$$\langle \cdot,\cdot \rangle
  \co H_2(\tilde{C}) \times H_2(\tilde{C},\partial\tilde{C} \cup \tilde{\nu})
  \to \Lambda$$
by
$$\langle x,y \rangle
  = \sum_{a,b \in \Z} (x \cdot q^a t^b y)q^a t^b.$$
A similar definition gives a pairing
$$\langle \cdot,\cdot \rangle'
  \co H_2(\tilde{C},\tilde{\nu}) \times H_2(\tilde{C},\partial\tilde{C})
  \to \Lambda.$$
To prove that these are well defined
requires some elementary homology theory.
See \cite[Appendix E]{aK96}.

For $x \in H_2(\tilde{C})$,
$y \in H_2(\tilde{C},\partial\tilde{C} \cup \tilde{\nu})$,
$\sigma \in B_n$, and $\lambda \in \Lambda$, we have
$$\langle \sigma x,\sigma y \rangle = \langle x,y \rangle,$$
and
$$\langle \lambda x,y \rangle = \lambda \langle x,y \rangle
                              = \langle x,\bar{\lambda}y \rangle,$$
where $\bar{\lambda}$ is the image of $\lambda$
under the automorphism of $\Lambda$
taking $q$ to $q^{-1}$ and $t$ to $t^{-1}$.
Similar identities hold for $\langle \cdot,\cdot \rangle'$.



\subsection{The Temperley-Lieb algebra}
\label{subsec:tl}

A good introduction to the Temperley-Lieb algebra
and its representation theory is \cite{bW95}.
We review some of the basics here.

Let $R$ be a domain containing an invertible element $q$
with a square root $q^{\half}$.
The Temperley-Lieb algebra $\TL_n(R)$
is the $R$-algebra given by generators $1,e_1,\dots,e_{n-1}$
and relations
\begin{itemize}
\item $e_i e_j = e_j e_i$ if $|i-j| > 1$,
\item $e_i e_j e_i = e_i$ if $|i-j| = 1$,
\item $e_i^2 = (-q^{\half}-q^{-\half})e_i$.
\end{itemize}
We simply write $\TL_n$ if $R$ is understood.
There is a map from $B_n$ to $\TL_n$
given by $\sigma_i \mapsto 1+q^{\half} e_i$.
Thus any representation of $\TL_n$
gives rise to a representation of $B_n$.

\begin{figure}
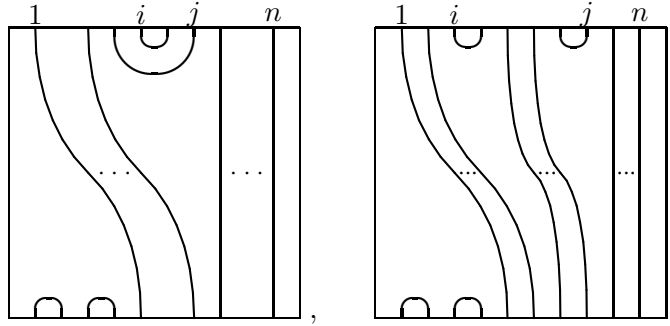

\centering
\sijclose,\qquad \sijfar
\caption{$s_{i,j}$ for $j=i+2$ and $j>i+2$.}
\label{fig:sij}
\end{figure}

Let $M$ be the left-ideal of $\TL_n$ generated by $e_1 e_3$,
and let $N$ be the left-ideal generated by $\{e_5,\dots,e_{n-1}\}$.
The representation of $\TL_n$ corresponding to the partition $(n-2,2)$
is the $\TL_n$-module $S = M/(M \cap N)$.
This is called $S(n,2)$ in \cite{bW95}.
It is a free $R$-module with basis
$$s_{i,j} = (e_i \dots e_3 e_2)(e_{j-1} \dots e_5 e_4)(e_1 e_3)$$
for $1 \le i < j \le n$ such that $j > \max(i+1,3)$.
There is a diagrammatic interpretation of $\TL_n$,
in which $s_{i,j}$ corresponds to the diagram shown in
Figure~\ref{fig:sij}.

\section{Some surfaces}
\label{sec:surfaces}

We describe some immersed surfaces in $\tilde{C}$
representing elements of homology and relative homology.
In each case we will describe a map from a surface to $C$
which can be lifted to a map to $\tilde{C}$.
We will not specify a choice of lift,
and we will not pay much attention to issues of orientation.
Thus the resulting element of homology or relative homology
will only be defined up to multiplication by a unit in $\Lambda$.
This will be sufficient for our purposes.

\subsection{Squares and triangles}

We describe some properly embedded surfaces representing elements of
relative homology modules such as $H_2(\tilde{C},\tilde{\nu})$.
These will be represented by one or two embedded edges in the disk.

Suppose $\alpha_1,\alpha_2 \co I \to D$ 
are disjoint embeddings of the interior of $I$ into $D_n$,
and map the endpoints of $I$ to puncture points 
(not necessarily injectively).
Let $f$ be the map from the interior of $I \times I$ to $C$
given by $f(x,y) = \{\alpha_1(x),\alpha_2(y)\}$.
A lift of $f$ to $\tilde{C}$
will represent an element of $H_2(\tilde{C},\tilde{\nu})$.
Similarly, we can define an element of 
$H_2(\tilde{C},\partial\tilde{C})$
corresponding to a pair of disjoint edges in $D_n$
with endpoints on $\partial D$.
Finally, we can define an element of
$H_2(\tilde{C},\partial\tilde{C} \cup \tilde{\nu})$
corresponding to a pair of edges which have a mixture of endpoints
on $\partial D$ and on puncture points.
For all of these examples,
call the resulting element of relative homology
the {\em square} corresponding to the edges $\alpha_1$ and $\alpha_2$.

Now suppose 
$\alpha \co I \to D$ is an embedding of the interior of $I$ into $D_n$,
and maps the endpoints of $I$ to the puncture points.
Define a map
$$f \co \{(x,y) \in I \times I : 0<x<y<1\} \to C$$
by $f(x,y)=\{f(x),f(y)\}$.
A lift of $f$ to $\tilde{C}$ will represent an element of
$H_2(\tilde{C},\tilde{\nu})$.
Similarly, we can obtain an element of
$H_2(\tilde{C},\partial\tilde{C}\cup\tilde{\nu})$
if we allow one or both endpoints of the edge to be on $\partial D$.
For all of these examples,
call the resulting element of relative homology
the {\em triangle} corresponding to the edge $\alpha$.

\subsection{A genus one surface}

\begin{figure}
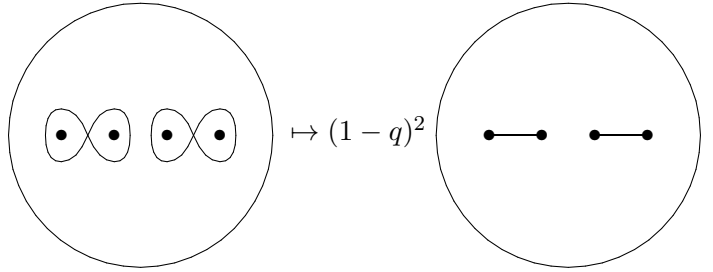

\centering
$\oooo \mapsto (1-q)^2 \oooosquare$.
\caption{A torus and a square.}
\label{fig:oooo_square}
\end{figure}

Suppose $\alpha_1,\alpha_2 \co S^1 \to D_n$ are disjoint figure-eights,
each going around two puncture points,
as in Figure~\ref{fig:oooo_square}.
Define a map $f$ from the torus $S^1 \times S^1$ to $C$ by
$f(x,y) = \{\alpha_1(x),\alpha_2(y)\}$.
Both the meridian and longitude of the torus
are mapped into the kernel of $\Phi$,
so $f$ lifts to a map 
$$\tilde{f} \co S^1 \times S^1 \to \tilde{C}.$$
We obtain an element $[\tilde{f}]$ of $H_2(\tilde{C})$.
The image of $[\tilde{f}]$ in $H_2(\tilde{C},\tilde{\nu})$
is $(1-q)^2$ times the square
as indicated in Figure~\ref{fig:oooo_square}.

\subsection{A genus two surface}

\begin{figure}
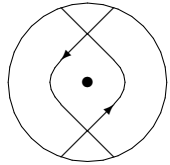

\centering
\ooodisk
\caption{The disk $B$ and edges $\alpha_1(I)$ and $\alpha_2(I)$.}
\label{fig:ooo_disk}
\end{figure}

\begin{figure}
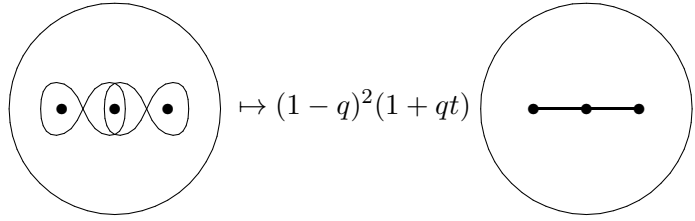

\centering
$\ooo \mapsto (1-q)^2(1+qt)\ooosquare$
\caption{A genus two surface and a square.}
\label{fig:ooo_square}
\end{figure}

Suppose $\alpha_1,\alpha_2 \co S^1 \to D_n$ are figure-eights
such that 
$\alpha_1$ passes around $p_i$ and $p_j$,
$\alpha_2$ passes around $p_j$ and $p_k$,
and $\alpha_1$ intersects $\alpha_2$ twice,
as in Figure~\ref{fig:ooo_square}.
Let $B \subset D$ be a disk centred at $p_j$,
containing the two points of intersection,
and meeting each of $\alpha_1$ and $\alpha_2$ in a single edge.
For $i=1,2$, let $I_i \subset S^1$ be the interval $\alpha_i^{-1}(B)$.
For convenience, assume $I_1 = I_2$,
and identify both with $I = [0,1]$,
oriented as shown in Figure~\ref{fig:ooo_disk}.

Let $T$ be the closure of
$$(S^1 \times S^1)\setminus (I_1 \times I_2).$$
Note that for $(x,y) \in T$ we have $\alpha_1(x) \neq \alpha_2(y)$.
We can therefore define $f \co T \to C$
by $f(x,y) = \{\alpha_1(x),\alpha_2(y)\}$.

Let $f_0$ be the restriction of $f$ 
to $\partial T = \partial(I \times I)$.
For $s \in I$,
let $f_s$ be the composition of $f_0$
with an anticlockwise rotation of $B$ 
by an angle of $s\pi$ about the centre $p_j$.
We can assume $B$ has a rotational symmetry so that 
$f_1(x,y) = f_0(y,x)$ for all $(x,y) \in \partial (I \times I)$.
Thus $f_0$ and $f_1$ represent the same loops,
but with opposite orientations.

We now build a closed genus two surface $\Sigma_2$ by gluing together
$\partial T \times I$ and two copies of $T$ as follows.
First glue $\partial T \times I$ to $T$ by
$((x,y),0) \sim (x,y)$.
Then glue $\partial T \times I$ to a second copy $T'$ of $T$ by
$((x,y),1) \sim (y,x)$.
Let $\Sigma_2$ be the surface so obtained.
Let $g \co \Sigma_2 \to C$ be given by
$g|_T = g|_{T'} = f$ and $g|_{(\partial T \times \{s\})} = f_s$.

The fundamental group of $\Sigma_2$
is generated by the meridian and longitude of $T$ and $T'$.
Each of these is mapped by $g$ into the kernel of $\Phi$,
so $g$ lifts to a map
$\tilde{g} \co \Sigma_2 \to \tilde{C}$.
We obtain an element $[\tilde{g}]$ of $H_2(\tilde{C})$.

We now compute the image of $[\tilde{g}]$
in $H_2(\tilde{C},\tilde{\nu})$.
For any $\epsilon>0$ we can assume that 
the disk $B$ has radius less than $\epsilon$,
and thus $g$ maps $\partial T \times I$ into $\nu_\epsilon$.
Then $\tilde{g}|_T$ represents an element of
$H_2(\tilde{C},\tilde{\nu}_\epsilon)$.
This element is $(1-q)^2$ times a square.
It remains to figure out how $\tilde{g}|_{T'}$ is related to $\tilde{g}|_T$.

Consider the path $\delta \co I \to \Sigma_2$ given by
$$\delta(s) = ((0,0),s) \in \partial T \times I.$$
This goes from $(0,0) \in T$ to $(0,0) \in T'$.
Now $g \circ \delta$ is a loop in $C$
in which the pair of points switch places by 
an anticlockwise rotation through an angle of $\pi$ around $p_j$.
Thus $\Phi(g \circ \delta) = qt$.
It follows that $\tilde{g}|_{T'} = qt \circ \tilde{g}|_T$.

Also note that $T'$ and $T$ inherit 
the same orientation from $\Sigma_2$.
This is because the ends of the annulus $\partial T \times I$
were attached with opposite orientations.
We conclude that the image of $[\tilde{g}]$ in $H_2(\tilde{C},\tilde{\nu})$ 
is $(1-q)^2(1+qt)$ times a square,
as shown in Figure~\ref{fig:ooo_square}. 

\subsection{A genus three surface}
\label{sec:genus_three}

\begin{figure}
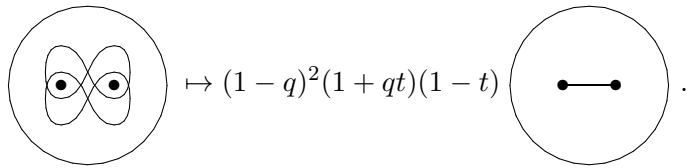

\centering
$\oo \mapsto (1-q)^2(1+qt)(1-t)\ootriangle$.
\caption{A genus three surface and a triangle.}
\label{fig:oo_triangle}
\end{figure}

Suppose $\alpha_1,\alpha_2 \co S^1 \to D_n$ are figure-eights,
both passing around $p_i$ and $p_j$,
and intersecting transversely at four points,
as in Figure~\ref{fig:oo_triangle}.
We construct a map from a genus three surface $\Sigma_3$ into $\tilde{C}$
by a slight modification of the procedure used above for $\Sigma_2$.
This time the surface $T$ will be a torus with two disks removed, one
for each of $p_i$ and $p_j$.
Two annuli are then needed to glue together two copies of $T$.
We obtain a genus three surface $\Sigma_3$ and a map $g\co\Sigma_3\to C$.
The image of the longitude, meridian, and both boundary components
of $T$ all lie in the kernel of $\Phi$,
so $g$ lifts to a map $\tilde{g} \co \Sigma_3\to\tilde{C}$.
This represents an element $[\tilde{g}]$ of $H_2(\tilde{C})$.

The image of $[\tilde{g}]$
in $H_2(\tilde{C},\tilde{\nu})$ is $(1+qt)(1-q)^2$ times a square,
which in turn is $(1+qt)(1-q)^2(1-t)$ times a triangle,
as shown in Figure~\ref{fig:oo_triangle}.

\section{A basis}
\label{sec:basis}

In this section, we give a new proof of the following theorem.

\begin{thm}
\label{thm:free}
$H_2(\tilde{C})$ is a free $\Lambda$-module of rank $\nchoosetwo$.
\end{thm}

This is originally due to Paoluzzi and Paris \cite{PP02}.
Our proof gives a more explicit description
of surfaces representing elements of a free basis for $H_2(\tilde{C})$.
We use the following basic result about $H_2(\tilde{C})$.

\begin{lem}
\label{lem:basic}
The vector space $\Q(q,t) \otimes H_2(\tilde{C})$
has dimension $\nchoosetwo$.
The natural map from $H_2(\tilde{C})$
to $\Q(q,t) \otimes H_2(\tilde{C})$ is injective.
\end{lem}

The proof uses a $2$-complex that is homotopy equivalent to $C$.
Various methods for constructing such a complex can be found in
\cite{rL90}, \cite{sB01}, \cite{PP02}, and \cite{rB02}.

\subsection{Proof of the theorem}

We now prove Theorem \ref{thm:free}.
For $1 \le i < j \le n$,
we define $v'_{i,j} \in H_2(\tilde{C},\tilde{\nu})$ as follows.
If $j-i > 2$,
let $v'_{i,j}$ be the square corresponding to the edges
$[p_i,p_{i+1}]$ and $[p_{j-1},p_j]$.
If $j-i = 2$ and $i > 1$,
let $v'_{i,j}$ be the square corresponding to
the edge $[p_i,p_{i+1}]$ 
and an edge from $p_{i-1}$ to $p_j$
whose interior lies in the lower half plane.
Let $v'_{1,3}$ be the square corresponding to the edges
$[p_1,p_2]$ and $[p_2,p_3]$.
If $j-i = 1$,
let $v'_{i,j}$ be the triangle corresponding to the edge $[p_i,p_j]$.

For $1 \le i < j \le n$,
let $v_{i,j}$ be an element of $H_2(\tilde{C})$
whose image in $H_2(\tilde{C},\tilde{\nu})$ is
\begin{itemize}
\item $(1-q)^2(1+qt)(1-t) v'_{i,j}$ if $j = i+1$,
\item $(1-q)^2(1+qt) v'_{i,j}$ if $i=1$ and $j=3$,
\item $(1-q)^2 v'_{i,j}$ otherwise.
\end{itemize}
Such an element exists by Section \ref{sec:surfaces}.
We will show that the $v_{i,j}$
form a free basis for $H_2(\tilde{C})$.

Every $x_{i,j} \in H_2(\tilde{C},\partial\tilde{C})$ 
be the square corresponding to a pair of vertical edges,
one passing just to the right of $p_i$,
the other passing just to the left of $p_j$,
and both having endpoints on $\partial D$.

\begin{lem}
\label{lem:pairings}
The following pairings are units in $\Lambda$.
\begin{itemize}
\item $\langle v'_{i,j},x_{i,j} \rangle'$ for $1 \le i < j \le n$,
\item $\langle v'_{i,i+2},x_{i-1,i+1} \rangle'$ for $i = 2,\dots,n-2$,
\item $\langle v'_{i,i+2},x_{i,i+1} \rangle'$ for $i = 2,\dots,n-2$.
\end{itemize}
All other pairings $\langle v'_{i',j'},x_{i,j} \rangle'$ are zero.
\end{lem}

\begin{proof}
Project the surfaces representing $v'_{i',j'}$ and $x_{i,j}$ to $C$.
If the resulting surfaces intersect transversely at precisely one point
then $\langle v'_{i',j'},x_{i,j} \rangle'$ is a unit in $\Lambda$.
If they do not intersect at all then the pairing is zero.
\end{proof}

For convenience, choose lifts and orientations so that
$$\langle v'_{i,j},x_{i,j} \rangle' = 1$$
for $1 \le i < j \le n$.

From Lemma \ref{lem:pairings} it is easy to show that
the $v'_{i,j}$ are linearly independent.
Hence the $v_{i,j}$ are linearly independent.
By Lemma \ref{lem:basic},
the $v_{i,j}$ span $\Q(q,t) \otimes H_2(\tilde{C})$ as a $\Q(q,t)$-module.
It remains to show that they span $H_2(\tilde{C})$ as a $\Lambda$-module.
In other words, we must prove the following.

\begin{lem}
\label{lem:lambda}
Let $c_{i,j} \in \Q(q,t)$ for $1 \le i < j \le n$ be such that
$$v = \sum_{1 \le i < j \le n} c_{i,j} v_{i,j}$$
lies in $H_2(\tilde{C})$.
Then $c_{i,j} \in \Lambda$ for all $1 \le i < j \le n$.
\end{lem}

We use the following facts about the $x_{i,j}$.

\begin{lem}
\label{lem:xright_multiple}
$x_{n-1,n}$ is a multiple of $(1-q)(1+qt)(1-t)$ 
in $H_2(\tilde{C},\partial\tilde{C}\cup\tilde{\nu})$.
\end{lem}

\begin{proof}
Let $\alpha \co I \to D$ be a straight edge from $\partial D$ to $p_n$.
Let $\alpha_1,\alpha_2 \co I \to D_n$ be disjoint edges 
that have endpoints on $\partial D$
and pass anticlockwise around $\alpha$,
such that $\alpha_1$ encloses $\alpha_2$.
Define $f \co I \times I \to C$ by $f(x,y) = \{\alpha_1(x),\alpha_2(y)\}$.
Then $f$ lifts to a map $\tilde{f}$ which represents $x_{2,3}$.

We can assume that for all $s \in I$ the four points
$\{\alpha_i(s),\alpha_i(1-s)\}$ lie within distance $\epsilon$ of each other,
where $\epsilon > 0$ is small.
Then $f$ maps the lines
$\{(x,x)\}$, $\{(x,-x)\}$, $\{(\half,y)\}$, $\{(x,\half)\}$
into $\nu_\epsilon$.
These lines cut $I \times I$ into eight pieces.
The restriction of $\tilde{f}$ to each of these pieces
represents an element of
$H_2(\tilde{C},\partial\tilde{C}\cup\tilde{\nu_\epsilon})$.
Each of these elements is a multiple of 
the triangle corresponding to $\alpha$,
as shown in Figure~\ref{fig:unionjack}.
\begin{figure}
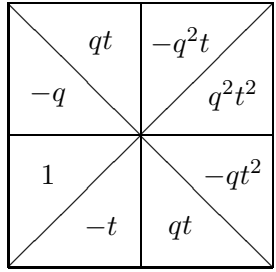

\centering
\unionjack
\caption{The function representing $x_{n-1,n}$ breaks into triangles.}
\label{fig:unionjack}
\end{figure}
Combining these, we see that $x_{n-1,n}$ is
$(1-q)(1+qt)(1-t)$ times the triangle corresponding to $\alpha$.
\end{proof}

\begin{lem}
\label{lem:x_multiple}
For all $1 \le i < j \le n$, $x_{i,j}$ is a multiple of $(1-q)^2$ 
in $H_2(\tilde{C},\partial\tilde{C}\cup\tilde{\nu})$.
\end{lem}

\begin{proof}
Use a similar argument to Lemma \ref{lem:xright_multiple},
as suggested by Figure~\ref{fig:mm}.
\begin{figure}
$$\mm = (1-q)(1-q^{-1})
              \left( \chopstickstwofour + q\chopsticksonefour + \dots \right)
$$
\caption{$x_{2,4}$ is a multiple of $(1-q)^2$.}
\label{fig:mm}
\end{figure}
\end{proof}

\begin{proof}[Proof of Lemma \ref{lem:lambda}]
First consider the case $n=3$. Then 
$$\langle v,x_{2,3} \rangle = (1-q)^2(1+qt)(1-t)c_{2,3}.$$
By Lemma \ref{lem:x_multiple}, $x_{2,3}$ is a multiple of $(1-q)^2$,
and hence of $(1-\bar{q})^2$.
Since $\langle \cdot,\cdot \rangle$ is sesquilinear,
it follows that $(1+qt)(1-t)c_{2,3} \in \Lambda$.
Similarly, Lemma \ref{lem:xright_multiple} implies that
$(1-q)c_{2,3} \in \Lambda$.
Since $\Lambda$ is a unique factorisation domain,
it follows that $c_{2,3} \in \Lambda$.
By a symmetrical argument, $c_{1,2} \in \Lambda$.
Now
$$v - (c_{1,2}v_{1,2} + c_{2,3}v_{2,3})$$
still lies in $H_2(\tilde{C})$,
so we can reduce to the case $v = c_{1,3}v_{1,3}$.
The following holds up to multiplication by a unit in $\Lambda$.
$$\langle v,\sigma_2 x_{2,3} \rangle = (1-t)(1-q)^2(1+qt)c_{1,3}.$$
Now $\sigma_2 x_{2,3}$ is a multiple of $(1-t)(1-q)(1+qt)$.
Thus $(1-q)c_{1,3}\in\Lambda$.
But
$$\langle v,x_{1,3} \rangle = (1-q)^2(1+qt)c_{1,3},$$
so $(1+qt)c_{1,3}\in\Lambda$.
Thus $c_{1,3} \in \Lambda$.

Now suppose $n > 3$.
For $i=1,\dots,n-2$ we have
$$\langle v,x_{i,n} \rangle = (1-q)^2c_{i,n}.$$
By Lemma \ref{lem:x_multiple}, $c_{i,n} \in \Lambda$.
Also
$$\langle v,x_{n-1,n} \rangle = (1-q)^2(1+qt)(1-t)c_{n-1,n}.$$
By Lemma \ref{lem:x_multiple}, $(1+qt)(1-t)c_{n-1,n} \in \Lambda$.
By Lemma \ref{lem:xright_multiple}, $(1-q)c_{n-1,n} \in \Lambda$.
Thus $c_{n-1,n} \in \Lambda$.
We can now subtract the terms $c_{i,n}v_{i,n}$,
and so reduce to the case $c_{i,n} = 0$ for all $i=1,\dots,n-1$.
Then $v$ represents an element of homology
of the preimage in $\tilde{C}$
of the space of unordered pairs of distinct points
in an $(n-1)$-times punctured disk.
The result now follows by induction on $n$.
\end{proof}

This completes the proof of Theorem \ref{thm:free}.

\subsection{The Krammer representation}
\label{sec:krammer}

There is some confusion as to the exact definition of the
``Lawrence-Krammer representation''.
In an attempt to clarify the situation,
we now compare and contrast a slightly different representation
which we will call the ``Krammer representation''.

Let the {\em Krammer representation}
be the following action of $B_n$
on a free $\Lambda$-module $V$ of rank $\nchoosetwo$
with basis $\{F_{i,j} : 1 \le i < j \le n\}$.
$$
\sigma_i(F_{j,k}) =
\left\{
\begin{array}{ll}
F_{j,k}
 &i \not \in \{j-1,j,k-1,k\}, \\
qF_{i,k}+(q^2-q)F_{i,j}+(1-q)F_{j,k}
 &i=j-1, \\
F_{j+1,k}
 &i=j \neq k-1, \\
qF_{j,i}+(1-q)F_{j,k}+(1-q)qtF_{i,k}
 &i=k-1\neq j, \\
F_{j,k+1}
 &i=k, \\
-tq^2 F_{j,k}
 &i=j=k-1.
\end{array}
\right.
$$
This action is given in \cite{dK00}, except that Krammer's $-t$ is my $t$.
It is also in \cite{sB01}, but with a sign error.
The name ``Krammer representation'' was chosen because
Krammer seems to have initially found this
independently of Lawrence and without any use of homology.

The vector spaces
$\Q(q,t) \otimes V$ and $\Q(q,t) \otimes H_2(\tilde{C})$
are isomorphic representations of $B_n$,
by \cite[Theorem 4.1]{sB01}.
The isomorphism is given by
$$F_{i,j} \mapsto 
(\sigma_{i-1} \dots \sigma_2 \sigma_1)
(\sigma_{j-1} \dots \sigma_3 \sigma_2) v_{1,2}.$$
However it is shown in \cite{PP02}
that the representations $V$ and $H_2(\tilde{C})$
are not isomorphic for $n \ge 3$.
The distinction becomes important when we specialise $q$ and $t$
to values which are not algebraically independent,
as we will in the next section.

It is possible to compute 
the matrices for the Lawrence-Krammer representation
with respect to the basis $\{v_{i,j}\}$.
We will not do this since they are quite complicated.
They must be conjugate to those of the Krammer representation
when considered as matrices over $\Q(q,t)$,
but not over $\Lambda$.

I conjecture that the Krammer representation
is the action of $B_n$ on $H_2(\tilde{C},\tilde{\nu})$.
This can be reduced to showing that
$\Q(q,t) \otimes H_2(\tilde{C},\tilde{\nu})$ has dimension $\nchoosetwo$.

\section{A representation of the Temperley-Lieb algebra}
\label{sec:tl}

The aim of this section is to prove the following.

\begin{thm}
\label{thm:tl}
Let $R$ be a domain containing invertible elements $q$ and $t$,
and let
$$\iota \co R \otimes H_2(\tilde{C})
        \to R \otimes H_2(\tilde{C},\tilde{\nu})$$
be induced by the natural map from homology to relative homology.
If $1+qt = 0$, $q$ has a square root, $q^2 \neq 1$, and $q^3 \neq 1$
then the image of $\iota$ is the representation of $\TL_n(R)$
corresponding to the partition $(n-2,2)$.
\end{thm}

The main example to keep in mind is
$R = \Z[q^{\pm \half}]$ with $t=-q^{-1}$.

In \cite{rL90}, Lawrence constructed the representation
corresponding to the partition $(n-k,k)$
for $k=1,\dots,\lfloor n/2 \rfloor$.
This was generalised still further in \cite{rL96}
to give the Hecke algebra representation
corresponding to any partition of $n$.
Theorem~\ref{thm:tl} only covers the case $(n-2,2)$,
but also has a number of advantages over the work of Lawrence.
Firstly, it gives a more elementary description
of the required quotient of $R \otimes H_2(\tilde{C})$.
Secondly, it works over a fairly general ring,
whereas Lawrence worked over $\C$
and used the matrices we gave in Section~\ref{sec:krammer}
for the ``Krammer representation''.
Finally, our proof gives an explicit isomorphism
between the two representations.
It is to be hoped that these advantages
can also be generalised to an arbitrary partition of $n$.

\subsection{A basis}

For $1 \le i < j \le n$,
let $v_{i,j}$, $v'_{i,j}$ and $x_{i,j}$
be as defined in Section~\ref{sec:basis}.
Then
$$\iota(v_{i,j}) =
  \left\{
    \begin{array}{ll}
    0 & \mbox{if}\, j=i+1 \,\mbox{or}\, j=3,\\
    (1-q)^2 v'_{i,j} & \mbox{otherwise}.
    \end{array}
  \right.
$$
The pairing $\langle \cdot,\cdot \rangle'$ can be extended to a map
$$\langle \cdot,\cdot \rangle'
  \co (R \otimes H_2(\tilde{C},\tilde{\nu}))
  \times H_2(\tilde{C},\partial\tilde{C})
  \to R.$$
Lemma \ref{lem:pairings} still holds,
so the $v'_{i,j}$ are linearly independent
in $R \otimes H_2(\tilde{C},\tilde{\nu})$.
Thus the image of $\iota$ is the free module with basis
$$\{(1-q)^2 v'_{i,j} : j > \max(i+1,3)\}.$$
Let $H$ be the free module with basis
$$\{v'_{i,j} : j > \max(i+1,3)\}.$$
To ease the notation, we will work with $H$ instead of the image of $\iota$,
since the two are isomorphic.

Let $K$ be the field of fractions of $R$.
Then $K \otimes H$ is a vector space of dimension $n(n-3)/2$,
and contains $H$ as an embedded submodule.
It is sometimes easier to prove that a relation holds in $H$
by proving that it holds in $K \otimes H$.

\subsection{The Hecke algebra}

We now prove that $H$ is a representation of the Hecke algebra.

\begin{defn}
The Hecke algebra $\cH_n(R)$, or simply $\cH_n$,
is the $R$-algebra given by generators $1,\sigma_1,\dots,\sigma_{n-1}$
and relations
\begin{itemize}
\item $\sigma_i \sigma_j = \sigma_j \sigma_i$ if $|i-j| > 1$,
\item $\sigma_i \sigma_j \sigma_i = \sigma_i \sigma_j \sigma_i$ 
      if $|i-j| = 1$,
\item $(\sigma_i - 1)(\sigma_i + q) = 0$.
\end{itemize}
\end{defn}

Thus $\cH_n$ is the group algebra $RB_n$ of the braid group
modulo the relations
$$(\sigma_i - 1)(\sigma_i + q) = 0.$$

\begin{lem}
\label{lem:hecke}
$(\sigma_i-1)(\sigma_i+q)$ acts as the zero map on $H$.
\end{lem}

\begin{proof}
Since the $\sigma_i$ are all conjugate,
we need only check the case $i = n-1$.
To do this, we find a basis for $K \otimes H$
consisting of eigenvectors of $\sigma_{n-1}$,
each having eigenvalue $1$ or $-q$.

The elements $v'_{i,j}$ for $\max(i+1,3) < j \le n-2$
are linearly independent eigenvectors of $\sigma_{n-1}$
with eigenvalue $1$.
There are $(n-2)(n-5)/2$ of these.

Let $\alpha$ be a simple closed curve
based at $p_{n-2}$ enclosing $p_{n-1}$ and $p_n$.
For $i=1,\dots,n-4$,
let $u_i$ be the square corresponding to
the edges $[p_i,p_{i+1}]$ and $\alpha$.
Let $u_{n-3}$ be the square corresponding to
the edge $[p_{n-3},p_{n-2}]$
and a simple closed curve
based at $p_{n-4}$ enclosing $p_{n-3},\dots,p_n$.
These are eigenvalues of $\sigma_{n-1}$ with eigenvalue $1$.
For $i,i' = 1,\dots,n-3$,
the following identities hold
up to multiplication by a unit in $\Lambda$.
\begin{itemize}
\item $\langle u_i,x_{i,n} \rangle' = 1-q$,
\item $\langle u_{i},x_{i',n} \rangle' = 0$ for $i' \neq i$,
\item $\langle v'_{i,j},x_{i',n} \rangle' = 0$
      for $1 \le i < j \le n-2$.
\end{itemize}
Thus the $u_i$ are linearly independent
and not in the span of $\{ v'_{i,j} : j \le n-2\}$.
Now let $u$ be the square corresponding to
$\alpha$ and a simple closed curve based at $p_{n-3}$ enclosing $\alpha$.
The following identities hold up to multiplication by a unit in $\Lambda$.
\begin{itemize}
\item $\langle u,x_{n-2,n-1} \rangle' = (1-q^2)(1+q^2t)(1-t)$,
\item $\langle u_i,x_{n-2,n-1} \rangle' = 0$ for $i = 1,\dots,n-3$,
\item $\langle v'_{i,j},x_{n-2,n-1} \rangle' = 0$ for $j \le n-2$.
\end{itemize}
We conclude that
$$\cB_1' = \{v_{i,j} : \max(i+1,3) < j \le n-2\}
           \cup \{u_1,\dots,u_{n-3}\}
           \cup \{u\}$$
is a linearly independent set of eigenvectors of $\sigma_{n-1}$
with eigenvalue $1$.
Each $v \in \cB_1$ lies in $H$,
as can be seen by constructing a surface in $\tilde{C}$
whose image in $H_2(\tilde{C},\tilde{\nu})$ is $(1-q)^2v$. 

For $i=1,\dots,n-3$, note that
$v'_{i,n}$ must be an eigenvector of $\sigma_{n-1}$,
because the function $\sigma_{n-1}$ can be chosen to fix setwise
the square in $C$ whose lift represents $v'_{i,n}$.
To find the eigenvalue,
note that 
$$\sigma_{n-1}^{-1}x_{i,n} = x_{i,n-1} - q^{-1}x_{i,n},$$
as shown in Figure~\ref{fig:eigenvalue},
\begin{figure}
$$\xonen = x_{1,n-1} - q^{-1} x_{1,n}$$
\caption{$\sigma_{n-1}^{-1}x_{i,n} = x_{i,n-1} - q^{-1}x_{i,n}$.}
\label{fig:eigenvalue}
\end{figure}
so
$$\langle \sigma_{n-1}v'_{i,n},x_{i,n} \rangle
  = \langle v'_{i,n},\sigma_{n-1}^{-1} x_{i,n} \rangle 
  = -q \langle v'_{i,n},x_{i,n} \rangle.
$$
Thus
$$\cB_{-q} = \{v'_{i,n} : i=1,\dots,n-3 \}.$$
is a set of eigenvectors of $\sigma_{n-1}$ with eigenvalue $-q$.
They are linearly independent, and since $-q \neq 1$,
they are not in the span of $\cB_1$.
By a dimension count,
$\cB_1 \cup \cB_{-q}$ is a basis for $K \otimes H$. 
Each element of this basis is annihilated by
$(\sigma_{n-1}-1)(\sigma_{n-1}+q)$,
so we are done.
\end{proof}

\subsection{The Temperley-Lieb algebra}

We now prove that $H$ is a representation of 
the Temperley-Lieb algebra $\TL_n$.
Consider the map from $\cH_n$ to $\TL_n$
given by $\sigma_i \mapsto 1+q^{\half}e_i$.
Using the presentation of the $\TL_n$
given in Section \ref{subsec:tl},
one can check that this is well-defined
and has kernel generated by
$$z_{i,j}
  = \sigma_i \sigma_j \sigma_i - \sigma_i \sigma_j - \sigma_j \sigma_i
  + \sigma_i + \sigma_j - 1$$
for $|i-j| = 1$.

\begin{lem}
\label{lem:tl}
$z_{i,j}$ acts as the zero map on $H$.
\end{lem}

\begin{proof}
The $z_{i,j}$ are all conjugate to each other,
so it suffices to prove that $z_{n-2,n-1}$
acts as the zero map on $H$.
We have
\begin{eqnarray*}
z_{n-2,n-1}&=&(\sigma_{n-2}\sigma_{n-1} - \sigma_{n-1} + 1)(\sigma_{n-2}-1)\\
           &=&(\sigma_{n-1}\sigma_{n-2} - \sigma_{n-2} + 1)(\sigma_{n-1}-1).
\end{eqnarray*}
Thus it suffices to find a basis for $K \otimes H$,
each of whose elements is
an eigenvector of either $\sigma_{n-1}$ or $\sigma_{n-2}$
with eigenvalue $1$.

Let $\cB_1$ be
the linearly independent set of eigenvectors of $\sigma_{n-1}$
with eigenvalue $1$
as defined in the proof of Lemma \ref{lem:hecke}.

Let $\alpha'$ be a simple closed curve
based at $p_n$ and enclosing $p_{n-1}$ and $p_{n-2}$.
For $i=1,\dots,n-4$,
let $w_i$ be the square corresponding to the edges
$[p_i,p_{i+1}]$ and $\alpha'$.
Let $w_{n-3}$ be the square corresponding to $\alpha'$
and a simple closed curve based at $p_{n-3}$ and enclosing $\alpha'$.
The $w_i$ all lie in $H$,
and are eigenvectors for $\sigma_{n-2}$ with eigenvalue $1$.

For $i=1,\dots,n-3$,
let $\xi_i \co I \to D_n$ be a vertical edge
with endpoints on $\partial D$,
passing just to the right of $p_i$,
and let $\gamma \co S^1 \to D_n$ be a figure-eight
going around $p_{n-1}$ and $p_n$ 
in a small regular neighbourhood of $[p_{n-1},p_n]$.
See Figure~\ref{fig:yone}.
\begin{figure}
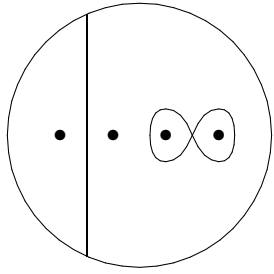

\centering
\yone
\caption{The arcs $\xi_1$ and $\gamma$ used to define $y_1$.}
\label{fig:yone}
\end{figure}
Define
$$f_i \co I \times S^1 \to C$$
by $f_i(x,y) = \{\xi_i(x),\gamma(y)\}$.
This lifts to $\tilde{C}$.
Let $y_i$ be the corresponding element of $H_2(\tilde{C},\partial\tilde{C})$.
For $i,j=1,\dots,n-3$,
the following holds up to multiplication by a unit in $\Lambda$.
$$
\langle w_i,y_i \rangle' =
\left\{
\begin{array}{ll}
1-q^3              & i=j\neq n-3 \\
(1-q^3)(1-q^3 t^2) & i=j=n-3 \\
0                  & i \neq j
\end{array}
\right.
$$
Further, $\langle v,y_i \rangle' = 0$ for all $v \in \cB_1$.
Thus
$$\{w_1,\dots,w_{n-3}\} \cup \cB_1$$
is linearly independent.
By a dimension count, it forms a basis of $K \otimes H$.
Each vector in this basis is annihilated by $z_{n-2,n-1}$,
so we are done.
\end{proof}

\subsection{The $(n-2,2)$ representation}

We have shown that $H$ can be considered as a representation of $\TL_n$.
It remains to show that it is isomorphic to the representation $S$
corresponding to the partition $(n-2,2)$.
Recall that $S$ is the $\TL_n$-module $M/(M \cap N)$, where
$M$ is the left-ideal of $\TL_n$ generated by $e_1 e_3$ and
$N$ is the left-ideal generated by $\{e_5, \dots, e_{n-1}\}$.

Let $\psi \co \TL_n \to H$ be the unique map such that $\psi(1) = v'_{1,4}$.
Then
$$\psi(e_3) = q^{-\half}(\sigma_3-1)(v'_{1,4})
           = (-q^{\half} - q^{-\half})v'_{1,4}.$$
Similarly $\psi(e_1) = (-q^{\half}-q^{-\half})v'_{1,4}$.
Thus 
$$\psi(e_1 e_3) = (q^{\half}+q^{-\half})^2 v'_{1,4}.$$
Let $\psi' \co M \to N$ be given by
$$\psi' = (q^{\half} + q^{-\half})^{-2} \psi|_M.$$

For $i=5,\dots,n-1$ we have that
$$\psi(e_i) = q^{-\half}(\sigma_i - 1)(v'_{1,4}) = 0.$$
Thus $\psi(N) = 0$,
so $\psi'(M \cap N) = 0$.
We therefore have an map $\phi \co S \to H$ induced by $\psi'$.

Recall that $S$ has a basis 
$$s_{i,j} = (e_i \dots e_2)(e_{j-1}\dots e_4)(e_1 e_3)$$
for $1 \le i < j \le n$ such that $j > \max(i+1,3)$.
We will show that $\phi(s_{i,j}) = v'_{i,j}$
up to multiplication by a unit in $\Lambda$.

\begin{lem}
For $j=5,\dots,n$ we have $(\sigma_{j-1} - 1)v'_{1,j-1} = v'_{1,j}$
up to multiplication by a unit in $\Lambda$.
\end{lem}

\begin{proof}
We can assume that the edge $\sigma_{j-1} ([p_{j-2},p_{j-1}])$
passes within $\epsilon$ of $p_{j-1}$.
Then we see that
$$\sigma_{j-1} v'_{1,j-1} = \lambda v_{1,j-1} + \mu v_{1,j}$$
for some units $\lambda,\mu \in \Lambda$.
To show that $\lambda = 1$,
note that $\sigma_{j-1} x_{1,j-1} = x_{1,j-1}$, so
\begin{eqnarray*}
\langle \sigma_{j-1} v'_{1,j-1},x_{1,j-1} \rangle'
   &=& \langle \sigma_{j-1} v'_{1,j-1},\sigma_{j-1} x_{1,j-1} \rangle' \\
   &=& \langle v'_{1,j-1},x_{1,j-1} \rangle'.
\end{eqnarray*}
\end{proof}

\begin{lem}
For $i=2,\dots,j-1$ we have $(\sigma_i - 1)v'_{i-1,j} = v'_{i,j}$
up to multiplication by a unit in $\Lambda$.
\end{lem}

\begin{proof}
In the case $i < j-1$,
the proof is almost identical to that of the previous lemma.
Suppose $i = j-1$.
Since everything takes place in a four-times punctured disk,
we can assume $j = n = 4$.
Then $H$ has a basis $\{v'_{1,4},v'_{2,4}\}$, so
$$\sigma_2 v'_{1,4} = \lambda v'_{1,4} + \mu v'_{2,4}$$
for some $\lambda,\mu \in \Lambda$.
To see that $\lambda = 1$,
use the same argument as in the proof of the previous lemma.
To see that $\mu$ is a unit in $\Lambda$,
note that it is equal to $\langle \sigma_2 v'_{1,4}, x_{2,4} \rangle'$,
and the surfaces representing $\sigma_2v'_{1,4}$ and $x_{2,4}$
intersect transversely exactly once in $C$.
\end{proof}

By these two lemmas we have $\phi(s_{i,j}) = v'_{i,j}$
up to multiplication by a unit in $\Lambda$.
This implies that $\phi$ is an isomorphism,
which completes the proof of Theorem \ref{thm:tl}.

\providecommand{\bysame}{\leavevmode\hbox to3em{\hrulefill}\thinspace}
\providecommand{\MR}{\relax\ifhmode\unskip\space\fi MR }
\providecommand{\MRhref}[2]{%
  \href{http://www.ams.org/mathscinet-getitem?mr=#1}{#2}
}
\providecommand{\href}[2]{#2}

\end{document}